\documentclass[12pt,reqno]{amsart}
\usepackage{amsfonts}
\usepackage{color}
\usepackage{epsfig}
\usepackage{fancyhdr}
\usepackage[latin1]{inputenc}
\usepackage{amsmath}
 \usepackage{mathrsfs}
\usepackage{thmdefs}
\usepackage{natbib}
\usepackage{graphicx}
\usepackage[colorlinks,citecolor=blue,urlcolor=blue]{hyperref}
\usepackage{amssymb}
\usepackage{dsfont}
\usepackage{fancyhdr}
\usepackage{graphicx}

\theoremstyle{definition} \theoremstyle{remark}
\numberwithin{equation}{section}

\renewcommand{\cite}{\citet}
\numberwithin{equation}{section}

\begin{document}

\author{Salim BOUZEBDA$^{*}$ and   Amor KEZIOU$^{**}$}

\address{$^{*}$LSTA-Universit\'e Paris 6, 175, rue du Chevaleret,
Bo\^ite 158, 75013 Paris, France.\\
 salim.bouzebda@upmc.fr}

\address{$^{**}$Laboratoire de Math\'ematiques (FRE 3111) CNRS,
Universit\'e de Reims and LSTA-Universit\'e Paris 6. UFR Sciences,
Moulin de la Housse, B.P. 1039,
                   51687 Reims,
                   France.
amor.keziou@upmc.fr
  }

\date{February  2009.}

\title{New estimates and tests of
independence in semiparametric copula models}

\maketitle

\begin{abstract}{\it  We introduce new estimates and tests of independence in
copula models with  unknown margins using $\phi$-divergences and
the duality technique. The asymptotic laws of the estimates and
the test statistics are established both when the parameter is an
interior or a boundary value of the parameter space. Simulation
results show that the choice of $\chi^2$-divergence has good
properties in terms of efficiency-robustness.}\\

\end{abstract}

\section{Introduction and motivations}
\noindent Copulas are a useful tool to model dependent data as
they allow to separate the dependence properties of the data from
their marginal properties and to construct multivariate models
with marginal distributions of arbitrary form. In particular,
parametric models for copulas with unknown margins have been
intensively investigated during the last decades. In the
monographs by \cite{Nelsen1999} and \cite{Joe1997} the reader
finds detailed accounts of the theory as well as surveys of
commonly used copulas.
\vskip7pt \noindent
It is known that some commonly used
dependence measures such as Pearson's correlation coefficient,
Kendall's tau and Spearman's rho cannot completely capture the
dependence structure among variables. Copulas have become
 popular in applied statistics, because
of the fact that they constitute a flexible and robust way to
model dependence between the margins of random vectors. \vskip7pt
\noindent In this framework, semiparametric inference methods,
based on \emph{pseudo-likelihood}, have been applied to copulas by
a number of authors (see, e.g., \cite{Shih_Louis1995},
\cite{Wang_Ding2000}, \cite{tsukahara2005} and the references
therein). Throughout the available literature, investigations on
the asymptotic properties of parametric estimators, as well as the
relevant test statistics, have privileged the case where the
parameter is an interior point of the admissible domain. However,
for most parametric copula models of interest, the boundaries of
the admissible parameter spaces include some important parameter
values, typically among which, that corresponding to the
independence of margins. We find in \cite{Joe1997} many examples
of parametric copulas, for which marginal independence is verified
for some specific values of the parameter $\theta$, on the
boundary $\partial \mathbf{\Theta}$ of the admissible parameter
set $\mathbf{\Theta} \subseteq \mathbf{R}^p$, $p\geq 1$. \vskip7pt
\noindent This paper concentrates  on this specific problem. We
aim, namely, to investigate parametric inference procedures, in
the case where the parameter belongs to the boundary of the
admissible domain. In particular, it will become clear, that the
usual limit laws both for parametric copula estimators and test
statistics become invalid under these limiting cases, and, in
particular, under marginal independence. Motivated by this
observation, we will introduce a new semiparametric inference
procedure based on \emph{$\phi$-divergences} and the
\emph{duality} technique extending
 the paper by \cite{bouzebda-keziou2009b} to the general
context of $\phi$-divergences for multivariate copulas with
multivariate parameter. The proposed method extends the
pseudo-maximum likelihood procedure introduced by
\cite{Genest_Ghoudi_Rivest1995}. It will be seen that the last
method corresponds to the particular choice of the
$KL_m$-divergence. We obtain a class of estimates and test
statistics depending upon the divergence. We are interested by
comparing the proposed estimates (including the pseudo-maximum
likelihood one) in terms of efficiency and robustness according to
the choice of the divergence. We will show that the proposed
estimators, under suitable conditions, remain asymptotically
normal, even under the marginal independence assumption for
appropriate choice of the divergence. This will allow us to
introduce test statistics of independence, whose study will be
made, both under the null and the alternative hypotheses. Let
$$\mathbf{F}(x_{1},\ldots,x_{d}):=\mathbf{P}\{X_1\leq
x_1,\ldots,X_d\leq x_d\}$$ be a $d$-dimensional distribution
function, and $F_{i}(x_i) := P(X_{i} \leq x_i),~i=1,\ldots,d$, the
marginal distributions of $\mathbf{F}(\cdot)$. It is well known
since the work of  \cite{Sklar1959} that there exists a
distribution  function $\mathbf{C}(\cdot)$ on $[0,1]^d$ with
uniform marginals such that
\begin{gather}\label{skl} \mathbf{C}(\mathbf{u}) :=\mathbf{C}(u_{1},\ldots,u_{d}) :=
\mathbf{P}\left\{F_{1}(X_{1}) \leq u_{1},\ldots,F_{d}(X_{d}) \leq u_{d}\right\}.
\end{gather} See also \cite{deheuvels1979b,Deheuvels1981b,deheuvels1981,Deheuvelsref2009}, \cite{Moore1975},
\cite{ruchendorf2009} and \cite{Schweizer1991}. We can refer to \cite{Sklar1973}, where the
author sketches the proof of (\ref{skl}), develops some of its
consequences, and surveys some of the work on copulas. Formally,
copulas can be defined in the common way as follows.
\begin{definition}
An $d$-dimensional copula is a function $\mathbf{C}:[0,1]^d\rightarrow
[0,1]$ with the following properties
\begin{enumerate}
\item $\mathbf{C}(\cdot)$ is grounded, i.e., for every
$\mathbf{u}=(u_1,\ldots,u_d)$, $\mathbf{C}(\mathbf{u})=0$ if at least one
coordinate $u_i=0,~i=1,\ldots,d$;
\item $\mathbf{C}(\cdot)$, is $d$-increasing,
i.e., for every $\mathbf{u}\in[0,1]^d$ and $\mathbf{v}\in[0,1]^d$
such that $\mathbf{u}\leq \mathbf{v}$, the $\mathbf{C}$-volume
$V_\mathbf{C}[\mathbf{u}, \mathbf{v}]$ of the box $[\mathbf{u},
\mathbf{v}]$ is non negative (see \cite{Nelsen1999}); \item
$\mathbf{C}(1,\ldots,1,u_i,1,\ldots,1)=u_i$ for all $u_i\in[0,1]^d$,
$u_i=0,~i=1,\ldots,d$.
\end{enumerate}
\end{definition}
\noindent Many useful multivariate models for dependence between
$X_1,\ldots,X_d$ turn out to be generated by \emph{parametric}
families of copulas of the form
$\left\{\mathbf{C}_\theta;~\theta\in\mathbf{\Theta}\right\}$, typically indexed by a
vector valued parameter $\theta\in\mathbf{\Theta}\subseteq\mathbf{R}^p$
(see, e.g., \cite{Kimeldorf1075_1}, \cite{Kimeldorf1075_2},
\cite{Nelsen1999},
 and \cite{Joe1993} among others). In the sequel, we assume that
 $\mathbf{C}_\theta(\cdot)$ admits a density $\mathbf{c}_\theta(\cdot)$
 with respect to the Lebesgue measure $\lambda$ on $\mathbf{R}^d$, i.e.,
 $\mathbf{c}_\theta(\cdot)=\frac{\partial^d  }{\partial u_1\ldots \partial u_d}\mathbf{C}_\theta(\cdot)$.
The nonparametric approach to copula estimation has been initiated
by \cite{deheuvels1979b}, who introduced and investigated the
\emph{empirical copula process}. In addition,
\cite{Deheuvels1981b,deheuvels1981,Deheuvelsref2009} described the limiting
behavior of this empirical process see, also
\cite{fermanianradulovicdragan2004} and the references therein.
The empirical copula process has been studied in full generality
in   \cite{Gaenssler1987} and \cite{stute1984}.
\vskip7pt \noindent In the present paper, we consider the estimation
and test problems for semiparametric copula models with unknown
general margins. Let $(X_{1k},\ldots,X_{dk})^\top$, for
$k=1,\ldots,n$, be a $d$-variate sample with distribution function
$\mathbf{F}_{\theta_\mathbf{T},F_1,\ldots,F_d}(\cdot,\ldots,\cdot)=\mathbf{C}_{\theta_\mathbf{T}}(F_1(\cdot),\ldots,F_d(\cdot))$
where $\theta_\mathbf{T}\in\mathbf{\Theta}$ is used to denote the
true unknown value of the parameter. In order to estimate
$\theta_\mathbf{T}$, some semiparametric estimation procedures,
based on the maximization, on the parameter space
$\mathbf{\Theta}$, of properly chosen \emph{pseudo-likelihood}
criterion, have been proposed by \cite{Genest_Ghoudi_Rivest1995},
\cite{Oakes1994}, \cite{Shih_Louis1995}, \cite{sellieng1996},
\cite{Wang_Ding2000} and \cite{tsukahara2005} among others. In
each of these papers, some asymptotic normality properties are
established for
$$\sqrt{n}\big(\tilde{\theta}-\theta_\mathbf{T}\big)$$ where
$\tilde{\theta}=\tilde{\theta}_n$ denotes a properly chosen
estimator of $\theta_\mathbf{T}$. This is achieved, provided that
$\theta_\mathbf{T}$ lies in the \emph{interior}, denoted by
$\mathring{\mathbf{\Theta}}$, of the parameter space
$\mathbf{\Theta}\subseteq \mathbf{R}^p$. On the other hand, the
case where $\theta_\mathbf{T}\in
\partial \mathbf{\Theta}:=\overline{\mathbf{\Theta}}-\mathring{\mathbf{\Theta}}$ is a
\emph{boundary value} of $\mathbf{\Theta}$, has not been studied
in a systematical way until present. Moreover, it turns out that,
for the above-mentioned estimators, the asymptotic normality of
$\sqrt{n}\big(\tilde{\theta}-\theta_\mathbf{T}\big),$ may fail to
hold for $\theta_\mathbf{T}\in
\partial \mathbf{\Theta}$; indeed, under some regularity conditions, when $\theta$ is univariate,  we can
prove that the limit law is the distribution of
$Z\mathds{1}_{(Z\geq 0)}$ where $Z$ is a centered normal variable,
and that the limit law of the  generalized pseudo-likelihood ratio
statistic is a mixture of chi-square laws with one degree of
freedom and Dirac measure  at zero; see
\cite{bouzebda-keziou2008}. Furthermore, when the parameter is
multivariate, the derivation of the limit distributions under the
null hypothesis of independence, becomes much more complex; see
\cite{selflieng1987}. Also, the limit distributions are not
standard which yields formidable numerical difficulties to
calculate the critical value of the test.\vskip7pt \noindent We cite below some
examples of parametric copulas, for which marginal independence is
verified for some specific values of the parameter $\theta$, on
the boundary $\partial \mathbf{\Theta}$ of the admissible
parameter set $\mathbf{\Theta}$.  We start with examples for which
$\theta$ varies within subsets of $\mathbf{R}$. Such is the case
for the extreme value copulas, namely
\begin{equation}
\mathbf{C}_{A}(u_1,u_2):=\exp \left\{\log u_1u_2A\left(\frac{\log u_1}{\log u_1u_2}\right)\right\},
\end{equation}
where $A(\cdot)$ is a convex function on $[0,1]$,  satisfying
\begin{enumerate}
\item[]   $A: [0, 1]\mapsto[1/2, 1]$ such that $\max(t, 1-t)\leq
A(t)\leq 1$ for all $0\leq t\leq 1$.
\end{enumerate}
For
\begin{equation}\label{gumbel_copula_f1}
A(t):=A_\theta(t)=(t^\theta+(1-t)^\theta)^{1/\theta} ; ~~\theta \in [1,\infty[
\end{equation}
 we have \cite{Gumbel1960} family of copulas, which is one of the most
popular model used to model bivariate extreme values. For
\begin{equation}\label{3a}
A_\theta(t)=1-(t^{-\theta}+(1-t)^{-\theta})^{-1/\theta}; ~~\theta\in [0, \infty[
\end{equation}
we obtain \cite{Galambos1975} family of copulas. Finally for
\begin{equation}\label{4}
A_\theta(t)=t\Phi\left(\theta^{-1}+\frac{1}{2}\theta\log\left(\frac{t}{1-t}\right)\right)+(1-t)\Phi\left(\theta^{-1}-
\frac{1}{2}\theta\log\left(\frac{t}{1-t}\right)\right),
\end{equation}
where $\theta \in[0,\infty[$ and $\Phi(\cdot)$ denoting the standard normal
$N(0,1)$ distribution function, we obtain the \cite{Husler_Reiss1989} family of copulas.
A useful family of copulas, due to \cite{Joe1993}, is given, for
$0<u_1,u_2<1$, by
\begin{eqnarray}\label{2}
\mathbf{C}_\theta(u_1,u_2):=1-\left[(1-u_1)^\theta+(1-u_2)^\theta-(1-u_1)^\theta
(1-u_2)^\theta\right]^{1/\theta}; ~~\theta\in
[1,\infty[.
\end{eqnarray}
The
Gumbel-Barnett copulas are given, for $0<u_1,u_2<1$, by
\begin{gather}\label{5}
\mathbf{C}_{\theta}(u_1,u_2) := u_1u_2\exp\left\{-(1-\theta)(\log u_1)(\log
u_2)\right\}; ~~ \theta \in [0,1].
\end{gather}
The Clayton copulas of positive dependence are such that, for $0<u_1,u_2<1$,
\begin{equation}\label{6}
\mathbf{C}_\theta(u_1,u_2)=\left(u_1^{-\theta}+u_2^{-\theta}-1\right)^{-1/\theta}; ~~\theta \in ]0,\infty[.
\end{equation}
Parametric families of copulas with parameter $\theta$ varying in
$\mathbf{R}^p$, for some $p\geq 2$, include the following
classical examples. Below, we set
$\theta=\big(\theta_1,\theta_2\big)^{\top}\in\mathbf{R}^2$.
\begin{equation}\label{7}
\mathbf{C}_{\theta}(u_1,u_2):=\left\{1+\left[(u_1^{-\theta_1}-1)^
{\theta_2}+(u_2^{-\theta_1}-1)^{\theta_2}\right]^{1/\theta_2}
\right\}^{-1/\theta_1}, ~\theta \in ]0,\infty[\times [1,\infty[;
\end{equation}
\begin{eqnarray}\label{8}
&&\mathbf{C}_{\theta}(u_1,u_2):=\exp\Big\{-\Big[{\theta_2}^{-1}\log
\Big(\exp\left(-{\theta_2}(\log
u_1)^{\theta_1}\right)\\
&&\qquad+\exp\left(-{\theta_2}(\log u_2)^{\theta_1}\right)
-1\Big)\Big]^{1/\theta_1}\Big\},~\theta \in [1,\infty[\times
]0,\infty[.\nonumber
\end{eqnarray}
For other examples of the kind, we refer to
\cite{Joe1997}.\vskip5pt \noindent For each of the above examples,
the independence case
$\mathbf{C}_{\theta_\mathbf{T}}(u_1,u_2)=u_1u_2$ (or $A(t)=1$)
occurs at the boundary of the parameter space $\mathbf{\Theta}$,
i.e., when $\theta_\mathbf{T}=1$ for the models
(\ref{gumbel_copula_f1}), (\ref{2}) and (\ref{5}),
$\theta_\mathbf{T}=0$ for the models (\ref{3a}), (\ref{4}) and
(\ref{6}), $\theta_\mathbf{T}=(0,1)^\top$ for the bivariate
parameter model (\ref{7}), and $\theta_\mathbf{T}=(1,0)^\top$ for
the bivariate parameter model (\ref{8}). In the sequel, we will
denote by $\theta_0$ the value of the parameter (when it exists),
corresponding to the independence of the marginals, i.e., the
value of the parameter for which we
have\begin{equation*}\mathbf{C}_{\theta_0}(\mathbf{u}):=\prod_{i=1}^du_i,~
\text{ for all } \mathbf{u}\in (0,1)^d.\end{equation*} Hence,
$\theta_0=1$ for the models (\ref{gumbel_copula_f1}), (\ref{2})
and (\ref{5}), $\theta_0=0$ for the models (\ref{3a}), (\ref{4})
and (\ref{6}), $\theta_0=(0,1)^\top$ for the model (\ref{7}), and
$\theta_0=(1,0)^\top$ for the model (\ref{8}). Note that for the
models (\ref{3a}), (\ref{4}), (\ref{6}), (\ref{7}) and (\ref{8}),
 $\mathbf{C}_{\theta_0}(u_1,u_2)=u_1u_2$ is naturally defined to be the limit of $\mathbf{C}_{\theta}(\cdot)$
 when $\theta$ tends to $\theta_0$ with values in $\mathbf{\Theta}$.
 Recall that  $\mathbf{c}_\theta(\cdot):=\frac{\partial^d}{\partial
u_1\ldots\partial u_d}\,\mathbf{C}_\theta(\cdot)$ is the density of $\mathbf{C}_\theta(\cdot)$ and
we define $\mathbf{c}_{\theta_0}(\cdot)$ to be the limit
 of $\mathbf{c}_{\theta}(\cdot)$
 when $\theta$ tends to $\theta_0$ with values in $\mathbf{\Theta}$. Hence, we can show that for all the above models
 $\mathbf{c}_{\theta_0}(u_1,u_2)=1$ for all $0<u_1,u_2<1$. \vskip5pt \noindent In contrast with the
preceding examples, where $\theta_0\in\partial\mathbf{\Theta}$ is a
boundary value of $\mathbf{\Theta}$, the case where $\theta_0$ is an
interior point of $\mathbf{\Theta}$ may, at times, occur, but is more
seldom. An example where $\theta_0\in\mathring{\mathbf{\Theta}}$ is given
by the Farlie-Gumbel-Morgenstern (FGM) copula, defined by
\begin{equation}\label{FGM}
\mathbf{C}_\theta (u_1,u_2):=u_1u_2+\theta u_1
u_2(1-u_1)(1-u_2),~\theta\in\mathbf{\Theta} :=[-1,1],
\end{equation}
and for which
$\theta_0=0\in\mathring{\mathbf{\Theta}}=]-1,1[$.\vskip5pt
\noindent In the present article, we will treat parametric
estimation of $\theta_\mathbf{T}$, and tests of the independence
assumption $\theta_\mathbf{T}=\theta_0$. We consider both the case
where $\theta_0\in\mathring{\mathbf{\Theta}}$ is an interior point
of $\mathbf{\Theta}$, and the case where $\theta_0\in \partial
\mathbf{\Theta}$ is a boundary value of $\mathbf{\Theta}$. Our
approach is novel in this setting and it will become clear later
on from our results, that the asymptotic normality of the estimate
based on $\phi$-divergences holds, even under the independence
assumption, when, either, $\theta_0$ is an interior, or a boundary
point of $\mathbf{\Theta}$, independently of the dimension of the
parameter space. The proposed test statistics of independence
using $\phi$-divergences are also studied, under the null
hypothesis $\mathscr{H}_0$ of independence, as well as under the
alternative hypothesis. The asymptotic distributions of the test
statistics under the alternative hypothesis are used to derive an
approximation to the power functions. An application of the
forthcoming results will allow us to evaluate the sample size
necessary to guarantee a pre-assigned power level, with respect to
a specified alternative.  To establish our results, we use similar
arguments as those developed by  \cite{tsukahara2005} in
connection with the instrumental statements on rank statistics
established by \cite{Ruymgaart_Shorack_Zwet1972} and
\cite{Ruschendorf1976} among others, combined with a new technique,
(based on the law of iterated logarithm given in Lemma
\ref{lemma1} below) to show both existence and consistency  of our
estimates and test statistics.
In \S\ref{simula}, we investigate the finite-sample performance of
the newly proposed estimators.
To avoid interrupting the flow of the
presentation, all mathematical developments are relegated to the appendix.
\section{A new inference procedure}
\noindent Recall that the $\phi$-divergence between a bounded
signed measure $\mathbf{Q}$, and a probability $\mathbf{P}$ on
$\mathscr{D}$, when $\mathbf{Q}$ is absolutely continuous with
respect to $\mathbf{P}$, is defined by
$$D_\phi(\mathbf{Q},\mathbf{P}):=\int_{\mathscr{D}}
\phi\left(\frac{d\mathbf{Q}}{d\mathbf{P}}(\bf x)\right)~d\mathbf{P(x)},$$
where $\phi$ is a convex function from $]-\infty,\infty[$ to
$[0,\infty]$ with $\phi(1)=0$. We will consider only
$\phi$-divergences for which the function $\phi$ is strictly
convex and satisfies: the domain of $\phi$, ${\rm dom}\phi:=\{x
\in \mathbf{R}: \phi(x)<\infty\}$ is an interval with end points
$a_{\phi}<1<b_{\phi}$, $\phi(a_\phi)=\lim_{x\downarrow
a_\phi}\phi(x)$ and $\phi(a_\phi)=\lim_{x\uparrow b_\phi}\phi(x)$.
The Kullback-Leibler, modified Kullback-Leibler, $\chi^2$,
modified $\chi^2$ and Hellinger  divergences are examples of
$\phi$-divergences; they are obtained respectively for
$\phi(x)=x\log x-x+1$, $\phi(x)=-\log x+x-1$,
$\phi(x)=\frac{1}{2}(x-1)^2$,
$\phi(x)=\frac{1}{2}\frac{(x-1)^2}{x}$ and
$\phi(x)=2(\sqrt{x}-1)^2$. 
We extend the definition of these
divergences on the whole space of all bounded signed  measures via
the extension of the definition of the corresponding $\phi$
functions on the whole real space $\mathbf{R}$ as follows: when
$\phi$ is not well defined on $\mathbf{R}_-$ or well defined but
not convex on $\mathbf{R}$, we set $\phi(x)=+\infty$ for all
$x<0$. Observe for the $\chi^2$-divergence, the corresponding
$\phi$ function is defined on whole $\mathbf{R}$ and strictly
convex. We refer to \cite{Liese-Vajda1987}  for a systematic
theory of divergences. We denote by $\phi^*$ the Fenchel-Legendre
transform of the convex function $\phi$, i.e., the function
defined by
$$t\in \mathbf{R}\mapsto \phi^*(t):=\sup_{x\in \mathbf{R}}\left\{tx-\phi(x)\right\}.$$
From \cite{Rockafellar1970}, Section 26, we can prove that it is
strictly convex, its domain is an interval $(a_\phi^*, b_\phi^*)$
with
$$a_\phi^*<0<b_\phi^*,\quad a_\phi^*=\lim_{x\to
-\infty}\frac{\phi(x)}{x}, \quad b_\phi^*=\lim_{x\to
+\infty}\frac{\phi(x)}{x},$$ and it satisfies $\phi^*(0)=0$,
$$\phi^*(a_\phi^*)=\lim_{t\downarrow a_\phi^*}\phi^*(t)\quad
\text{ and }\quad \phi^*(b_\phi^*)=\lim_{t\uparrow
b_\phi^*}\phi^*(t).$$ Furthermore, it holds that $\phi$ is the
Fenchel-Legendre transform of $\phi^*$. In the sequel,  for all
$\theta$, we denote by $D_\phi(\theta,\theta_\mathbf{T})$ the
$\phi$-divergences between $\mathbf{C}_{\theta}(\cdot)$ and
$\mathbf{C}_{\theta_\mathbf{T}}(\cdot)$, i.e.,
\begin{equation}
D_\phi(\theta,\theta_\mathbf{T}):=\int_{\mathbf{I}}
\phi\left(\frac{d\mathbf{C}_{\theta}}{d\mathbf{C}_{\theta_\mathbf{T}}}
(\mathbf{u})\right)~d\mathbf{C}_{\theta_\mathbf{T}}(\mathbf{u})=
\int_{\mathbf{I}}\phi\left(\frac{\mathbf{c}_{\theta}(\mathbf{u})}
{\mathbf{c}_{\theta_\mathbf{T}}(\mathbf{u})}\right)~d\mathbf{C}_{\theta_\mathbf{T}}(\mathbf{u}),
\end{equation}
where $\mathbf{I}=(0,1)^d.$ Denote $\mathbf{C}_n(\cdot)$ the empirical copula associated to the data, i.e.,
\begin{equation}\label{ec}
\mathbf{C}_n(\mathbf{u}):=\frac{1}{n}\sum_{k=1}^{n}\prod_{i=1}^d\mathds{1}_{\left\{F_{in}(X_{ik})\leq
u_1\right\}},~\mathbf{u}\in I,
\end{equation}
and $$F_{in}(t):=
\left\{\frac{n}{n+1}\right\}\frac{1}{n}\sum_{k=1}^n
\mathds{1}_{]-\infty,t]}(X_{ik})=\frac{1}{n+1}\sum_{k=1}^n
\mathds{1}_{]-\infty,t]}(X_{ik}),~i=1,\ldots,d,$$ where
$\mathds{1}_A$ stands for the indicator function of the event $A$.
The rescaling by the factor $n/(n+1)$, avoids difficulties arising
from potential unboundedness of $\mathbf{c}_{\theta}(\mathbf{u})$
when one of  $u_i$'s tends to $1$. Observe that the plug-in
estimate
$$\int_{\mathbf{I}}\phi\left(\frac{d\mathbf{C}_{\theta}}{d\mathbf{C}_{n}}(\mathbf{u})\right)~d\mathbf{C}_{n}(\mathbf{u})$$
of $D_\phi(\theta,\theta_\mathbf{T})$ is not well defined since
$\mathbf{C}_\theta(\cdot)$ is not absolutely continuous with
respect to $\mathbf{C}_n(\cdot)$. In order to avoid this
difficulty,  and to estimate the divergences
$D_\phi(\theta,\theta_\mathbf{T})$ for a given
$\theta\in\mathbf{\Theta}$ in particular for $\theta=\theta_0$, we
will make use of  the dual representation of $\phi$-divergences
obtained by \cite{BK2005111} Theorem 4.4 and \cite{Keziou2003}
Theorem 2.3. By this, when $\phi$ is differentiable, we readily
obtain that $D_\phi(\theta_0,\theta_\mathbf{T})$ can be rewritten
into
\begin{equation}\label{eqn 2 dual}
D_\phi(\theta_0,\theta_\mathbf{T}):=\sup_{f\in\mathscr{F}}\left\{\int_{\mathbf{I}}
f~d\mathbf{C}_{\theta_0}-\int_{\mathbf{I}}
\phi^*(f)~d\mathbf{C}_{\theta_\mathbf{T}}\right\},
\end{equation}
where $\mathscr{F}$ is an arbitrary class of measurable functions
fulfilling  the following two conditions:
\begin{equation*}
\forall
f\in\mathscr{F},~\int |f|~d\mathbf{C}_{\theta_0}<\infty
\end{equation*} and
\begin{equation*}
\phi'(d\mathbf{C}_{\theta_0}/d\mathbf{C}_{\theta_\mathbf{T}})=\phi'(\mathbf{c}_{\theta_0}/\mathbf{c}_{\theta_\mathbf{T}})\in\mathscr{F}.
\end{equation*}
\noindent Furthermore, the sup in the above display is unique and
is achieved at
$f=\phi'(\mathbf{c}_{\theta_0}/\mathbf{c}_{\theta_\mathbf{T}}).$
Note that for the specific value $\theta_0$, corresponding to the
independence, we have $\mathbf{c}_{\theta_0}(\mathbf{u})=1,
~\forall \mathbf{u}\in \mathbf{I}.$ So, by the above statement, taking the
class of functions
$$\mathscr{F}=\left\{\mathbf{u}\in \mathbf{I}\mapsto \phi'\left(1/\mathbf{c}_\theta(\mathbf{u})\right);~\theta\in\mathbf{\Theta}\right\},$$
we obtain the formula
\begin{equation}
D_\phi(\theta_0,\theta_\mathbf{T})=
\sup_{\theta\in\mathbf{\Theta}}\left\{\int_{\mathbf{I}}\phi'\left(\frac{\mathbf{c}_{\theta_0}}
{\mathbf{c}_\theta}\right)d\mathbf{C}_{\theta_0}(\mathbf{u})-
\int_{\mathbf{I}}\left[\frac{\mathbf{c}_{\theta_0}}{\mathbf{c}_\theta}\phi'\left(\frac{\mathbf{c}_{\theta_0}}{\mathbf{c}_\theta}\right)-
\phi\left(\frac{\mathbf{c}_{\theta_0}}
{\mathbf{c}_\theta}\right)\right]d\mathbf{C}_{\theta_\mathbf{T}}(\mathbf{u})\right\}\nonumber\\
\end{equation}
\begin{equation}\label{formula 1}
=\sup_{\theta\in\mathbf{\Theta}}\left\{\int_{\mathbf{I}}\phi'\left(\frac{1}{\mathbf{c}_\theta}\right)du_1\ldots
du_d-
\int_{\mathbf{I}}\left[\frac{1}{\mathbf{c}_\theta}\phi'\left(\frac{1}{\mathbf{c}_\theta}\right)-
\phi\left(\frac{1}{\mathbf{c}_\theta}\right)\right]
d\mathbf{C}_{\theta_\mathbf{T}}(\mathbf{u})\right\},
\end{equation} whenever
$$\int_{\mathbf{I}}\left|\phi'\left(1/\mathbf{c}_\theta\right)\right|~du_1\ldots
du_d <\infty ~~ ~\mbox{ for all }~~ \theta\in\mathbf{\Theta}.$$
Furthermore, the sup is unique and reached at
$\theta=\theta_\mathbf{T}$. Hence, the divergence
$D_\phi(\theta_0,\theta_\mathbf{T})$ and the parameter
$\theta_\mathbf{T}$ can be estimated respectively by
\begin{equation}\label{estimateur 1 div }
\sup_{\theta\in\mathbf{\Theta}}\left\{\int_{\mathbf{I}}\phi'\left(\frac{1}{\mathbf{c}_\theta}\right)du_1\ldots
du_d-
\int_{\mathbf{I}}\left[\frac{1}{\mathbf{c}_\theta}\phi'\left(\frac{1}{\mathbf{c}_\theta}\right)-
\phi\left(\frac{1}{\mathbf{c}_\theta}\right)\right]d\mathbf{C}_{n}(\mathbf{u})\right\}
\end{equation}
and
\begin{equation}\label{estimateur 1 param}\arg\sup_{\theta\in\mathbf{\Theta}}\left\{\int_{\mathbf{I}}
\phi'\left(\frac{1}{\mathbf{c}_\theta}\right)~du_1\ldots du_d-
\int_{\mathbf{I}}\left[\frac{1}{\mathbf{c}_\theta}\phi'\left(\frac{1}{\mathbf{c}_\theta}
\right)-\phi\left(\frac{1}
{\mathbf{c}_\theta}\right)\right]~d\mathbf{C}_{n}(\mathbf{u})\right\},
\end{equation}
in which $\mathbf{C}_{\theta_\mathbf{T}}(\cdot)$ is replaced by
$\mathbf{C}_n(\cdot)$. Note that this class of estimates contains
the maximum pseudo-likelihood (MPL) estimator proposed by
\cite{Oakes1994}; it is obtained for the $KL_m$-divergence taking
$\phi(x)=-\log(x)+x-1$. Under some regularity conditions, we can
prove that these estimates are consistent and asymptotically
normal in the same way as the MPL estimate when the parameter
$\theta_\mathbf{T}$ is an interior point of the parameter space
$\mathbf{\Theta}.$ The interest of divergence remains in the fact
that a properly choice of the divergence may ameliorate the MPL
estimator in terms of efficiency-robustness. The results in
\cite{bouzebda-keziou2008} show that, for
$\mathbf{\Theta}=[\theta_0,\infty)$, and when the true value
$\theta_\mathbf{T}$ of the parameter is equal to $\theta_0$
(corresponding to the independence assumption), the classical
asymptotic normality property of the MPL estimate is no more
satisfied. To overcome this difficulty,  in what follows, we
enlarge the parameter space $\mathbf{\Theta}$ into a wider space
$\mathbf{\Theta}_e\supset \mathbf{\Theta}$. This is tailored to
let $\theta_0$ become an interior point of $\mathbf{\Theta}_e$.
More precisely, set
\begin{equation}\label{def de Theta_e }
\mathbf{\Theta}_e :=\left\{\theta\in\mathbf{R}^p \text{ such that
} \int_{\mathbf{I}} \left|\phi'(1/\textbf{c}_\theta(\mathbf{u}))\right|~du_1\ldots
du_d <\infty\right\}.
\end{equation}
So, applying (\ref{eqn 2 dual}), with the class of functions
$$\mathscr{F}:=\left\{\mathbf{u}\in \mathbf{I}\mapsto
\phi'(1/{\textbf{c}}_\theta(\mathbf{u})); ~\theta\in\mathbf{\Theta}_e\right\},$$ we
obtain
\begin{equation}\label{formula 2}
D_\phi(\theta_0,\theta_\mathbf{T})=\sup_{\theta\in\mathbf{\Theta}_e}\left\{\int_{\mathbf{I}}\phi'\left(\frac{1}{\mathbf{c}_\theta}\right)d\mathbf{u}-
\int_{\mathbf{I}}\left[\frac{1}{\mathbf{c}_\theta}\phi'\left(\frac{1}{\mathbf{c}_\theta}\right)-
\phi\left(\frac{1}{\mathbf{c}_\theta}\right)\right]d\mathbf{C}_{\theta_\mathbf{T}}(\mathbf{u})\right\}.
\end{equation}
Furthermore, the sup in this display is unique and reached in
$\theta=\theta_\mathbf{T}.$ Hence, we propose to estimate
$D_\phi(\theta_0,\theta_\mathbf{T})$ by
\begin{equation}\label{estim div}
\widehat{D_\phi}(\theta_0,\theta_\mathbf{T}) := \sup_{\theta \in
\mathbf{\Theta}_e} \int_{\mathbf{I}}
\mathbf{m}(\theta,\mathbf{u})~d\mathbf{C}_{n}(\mathbf{u}),
\end{equation}
and to estimate the parameter $\theta_\mathbf{T}$ by
\begin{equation}\label{estimate}
\widehat\theta_n := \arg\sup_{\theta \in \mathbf{\Theta}_e}
\left\{\int_{\mathbf{I}}
\mathbf{m}(\theta,\mathbf{u})~d\mathbf{C}_{n}(\mathbf{u})\right\},
\end{equation}
where
\begin{eqnarray*}
\mathbf{m}(\theta,\mathbf{u})&:=&\int_{\mathbf{I}}\phi^{\prime}
\left(\frac{1}{\mathbf{c}_\theta (\mathbf{u})}\right)~d\mathbf{u}-
\left\{\phi^{\prime}\left(\frac{1}{\mathbf{c}_{\theta}(\mathbf{u})}\right)
 \frac{1}{\mathbf{c}_{\theta}(\mathbf{u})}-
 \phi\left(\frac{1}{\mathbf{c}_\theta(\mathbf{u})}\right)\right\}.
\end{eqnarray*}
In the sequel we denote by $\frac{\partial}{\partial
\theta}\mathbf{m}(\theta,\mathbf{u})$ the $p$-dimensional vector
with entries
$\frac{\partial}{\partial\theta_{i}}\mathbf{m}(\theta,\mathbf{u})$
and $\frac{\partial^2}{\partial
\theta^2}\mathbf{m}(\theta,\mathbf{u})$ the $p\times p$-matrix
with entries
$\frac{\partial^{2}}{\partial\theta_{i}\partial\theta_{j}}\mathbf{m}(\theta,\mathbf{u})$.
In what follows, we give some examples of divergences and the
associated estimates.

\subsection{Examples}
\begin{itemize}
  \item Our first example is the common used modified Kullback-Leibler divergence
  \begin{eqnarray*}
  \phi(x)&=&-\log x+x-1\\
\phi^\prime(x)&=&-\frac{1}{x}+1\\
x  \phi^\prime(x)-\phi(x)&=&\log x.
\end{eqnarray*}
The estimate of $D_{\rm KL_m}(\theta_0,\theta_\mathbf{T})$ is
given by
\begin{eqnarray*}
\widehat{D}_{\rm
KL_m}(\theta_0,\theta_\mathbf{T})&=&\sup_{\theta\in\mathbf{\Theta}_e}\left\{-\int_{\mathbf{I}}\log\left(\frac{1}{\mathbf{c}_\theta(\mathbf{u})}
\right)d\mathbf{C}_n(\mathbf{u})\right\}\\&=&\sup_{\theta\in\mathbf{\Theta}_e}\left\{\int_{\mathbf{I}}
\log\left(\mathbf{c}_\theta(\mathbf{u})\right)
d\mathbf{C}_n(\mathbf{u})\right\}
\end{eqnarray*}
and the estimate of the parameter $\theta_\mathbf{T}$ is given by
\begin{equation*}
\widehat\theta_n := \arg\sup_{\theta \in \mathbf{\Theta}_e}
\left\{\int_{\mathbf{I}}
\log\left(\mathbf{c}_\theta(\mathbf{u})\right)
d\mathbf{C}_n(\mathbf{u})\right\},
\end{equation*}
which coincides with the MPL one.
  \item The second one is the  Kullback-Leibler divergence
  \begin{eqnarray*}
  \phi(x)&=&x\log x-x+1\\
\phi^\prime(x)&=&\log x\\
x  \phi^\prime(x)-\phi(x)&=&x-1.
\end{eqnarray*}
The estimate of $D_{\rm KL}(\theta_0,\theta_\mathbf{T})$ is given
by
\begin{eqnarray*}
\widehat{D}_{\rm
KL}(\theta_0,\theta_\mathbf{T})&=&\sup_{\theta\in\mathbf{\Theta}_e}\left\{\int_{\mathbf{I}}\log\left(\frac{1}{\mathbf{c}_\theta(\mathbf{u})}
\right)d\mathbf{u}-\int_{\mathbf{I}}\left(\frac{1}{\mathbf{c}_\theta(\mathbf{u})}-1
\right)d\mathbf{C}_n(\mathbf{u})\right\}
\end{eqnarray*}
and the estimate of the parameter $\theta_\mathbf{T}$ is defined
as follows
\begin{equation*}
\widehat\theta_n := \arg\sup_{\theta \in \mathbf{\Theta}_e}
\left\{\int_{\mathbf{I}}\log\left(\frac{1}{\mathbf{c}_\theta(\mathbf{u})}
\right)d\mathbf{u}-\int_{\mathbf{I}}\left(\frac{1}{\mathbf{c}_\theta(\mathbf{u})}-1
\right)d\mathbf{C}_n(\mathbf{u})\right\}.
\end{equation*}
  \item The third one is the $\chi^2$-divergence
  \begin{eqnarray*}
  \phi(x)&=&\frac{1}{2}(x-1)^2\\
\phi^\prime(x)&=&x-1\\
x  \phi^\prime(x)-\phi(x)&=&\frac{1}{2}x-\frac{1}{2}.
\end{eqnarray*}
The estimate of $D_{\rm \chi^2}(\theta_0,\theta_\mathbf{T}) $ is given by
\begin{eqnarray*}
\widehat{D}_{\rm \chi^2}(\theta_0,\theta_\mathbf{T})&=&
\sup_{\theta\in\mathbf{\Theta}_e}\left\{\int_{\mathbf{I}}\left(\frac{1}{\mathbf{c}_\theta(\mathbf{u})}-1
\right)d\mathbf{u}\right.\\&-&\left.\int_{\mathbf{I}}\frac{1}{2}\left(\left(\frac{1}{\mathbf{c}_\theta(\mathbf{u})}\right)^2-1
\right)d\mathbf{C}_n(\mathbf{u})\right\}
\end{eqnarray*}
and the estimate of the parameter $\theta_\mathbf{T}$ is defined by
\begin{equation*}
\widehat\theta_n := \arg\sup_{\theta \in \mathbf{\Theta}_e}\left\{\int_{\mathbf{I}}\left(\frac{1}{\mathbf{c}_\theta(\mathbf{u})}-1
\right)d\mathbf{u}-\int_{\mathbf{I}}\frac{1}{2}\left(\left(\frac{1}{\mathbf{c}_\theta(\mathbf{u})}\right)^2-1
\right)d\mathbf{C}_n(\mathbf{u})\right\}.
\end{equation*}
\item The last example is the Hellinger divergence
\begin{eqnarray*}
  \phi(x)&=&2(\sqrt{x}-1)^2\\
\phi^\prime(x)&=&2-\frac{1}{\sqrt{x}}\\
x  \phi^\prime(x)-\phi(x)&=&2\sqrt{x}-2.
\end{eqnarray*}
The estimate of $D_{\rm H}(\theta_0,\theta_\mathbf{T})$ is given by
\begin{equation*}
\hspace{-.3cm}\widehat{D}_{\rm H}(\theta_0,\theta_\mathbf{T})=\sup_{\theta\in\mathbf{\Theta}_e}
\left\{\int_{\mathbf{I}}\left(2-2\sqrt{\mathbf{c}_\theta(\mathbf{u})}\right)d\mathbf{u}
-\int_{\mathbf{I}}2\left(\frac{1}{\sqrt{\mathbf{c}_\theta(\mathbf{u})}}-1
\right)d\mathbf{C}_n(\mathbf{u})\right\}
\end{equation*}
and the estimate of the parameter $\theta_\mathbf{T}$ is defined by
\begin{equation*}
\widehat\theta_n := \arg\sup_{\theta \in \mathbf{\Theta}_e}\left\{\int_{\mathbf{I}}\left(2-2\sqrt{\mathbf{c}_\theta(\mathbf{u})}\right)d\mathbf{u}
-\int_{\mathbf{I}}2\left(\frac{1}{\sqrt{\mathbf{c}_\theta(\mathbf{u})}}-1
\right)d\mathbf{C}_n(\mathbf{u})\right\}.
\end{equation*}
\end{itemize}
All the above examples are particular cases of the so-called
``\emph{power divergences}'', introduced by
\cite{Cressie-Read1984} (see also \cite{Liese-Vajda1987} Chapter
2), which are defined through the class of convex real valued
functions
\begin{equation*}\label{gamma convexfunctions}
x\in\mathbb{R}_{+}^*\rightarrow
\varphi_{\gamma}(x):=\frac{x^{\gamma }-\gamma x+\gamma -1}{\gamma
(\gamma -1)}
\end{equation*}
for $\gamma$ in $\mathbb{R}\backslash\left\{0,1\right\}$. The estimate of $D_{\gamma}(\theta_0,\theta_\mathbf{T})$ is given by
\begin{eqnarray*}
\widehat{D}_{\gamma}(\theta_0,\theta_\mathbf{T})&=&\sup_{\theta\in\mathbf{\Theta}_e}
\left\{\int_{\mathbf{I}}\frac{1}{\gamma-1}\left(\left(\frac{1}{\mathbf{c}_\theta(\mathbf{u})}\right)^{\gamma-1}-1
\right)d\mathbf{u}\right.\\&-&\left.\int_{\mathbf{I}}\frac{1}{\gamma}\left(\left(\frac{1}{\mathbf{c}_\theta(\mathbf{u})}\right)^\gamma-1
\right)d\mathbf{C}_n(\mathbf{u})\right\}
\end{eqnarray*}
and the parameter estimate  is defined by
\begin{eqnarray*}
\widehat\theta_n &:=& \arg\sup_{\theta \in \mathbf{\Theta}_e}
\left\{\int_{\mathbf{I}}\frac{1}{\gamma-1}\left(\left(\frac{1}{\mathbf{c}_\theta(\mathbf{u})}\right)^{\gamma-1}-1
\right)d\mathbf{u}\right.\\ &&-\left.\int_{\mathbf{I}}\frac{1}{\gamma}\left(\left(\frac{1}{\mathbf{c}_\theta(\mathbf{u})}\right)^{\gamma}-1
\right)d\mathbf{C}_n(\mathbf{u})\right\}.
\end{eqnarray*}

\begin{remark}Divergences measures have been intensively
used in estimation and test in the framework of the discrete
parametric models  with independent identically distributed data;
the estimates of the divergences and the parameter are obtained by
the plug-in method; see \cite{Liese-Vajda1987}  including the references
therein. For continuous parametric models the plug-in procedure
does not lead to well defined estimates; \cite{Keziou2003},
\cite{liesevajda2006}, \cite{Broniatowski-Keziou2008} introduce
new estimates and tests, using the dual representation of
divergences, extending the maximum likelihood procedure.
\end{remark}

\begin{remark}
We give an example of copulas for which the likelihood-based procedure fails.
We consider the Gumbel copulas $\mathbf{C}_\theta(\cdot)$ given in (\ref{gumbel_copula_f1}), it's corresponding density copula is defined by
    \begin{eqnarray}
    \mathbf{c}_\theta(u_1,u_2):=\mathbf{C}_\theta(u_1,u_2)(u_1u_2)^{-1}\frac{(\widetilde{u_1}\widetilde{u_2})^{(\theta-1)}}
    {(\widetilde{u_1}^\theta+\widetilde{u_2}^\theta)^{(2-1/\theta)}}
    \left[(\widetilde{u_1}^\theta+\widetilde{u_2}^\theta)^{(1/\theta)}+\theta-1\right],
    \end{eqnarray}
    where $\widetilde{x}=-\log x$.
We can show that $\mathbf{c}_\theta(\cdot)$ may takes negative
values for some $\theta\in \mathbf{\Theta}_e$. In fact
$\mathbf{c}_{0.7}(u_1,u_2)$ is negative for
$(u_1,u_2)\in[0.9,1]^2$, hence the likelihood function is not well
defined. The choice of the $\chi^2$-divergence is particularly
well adapted to this situation for example.
\end{remark}

\section{The asymptotic behavior of
the estimates}
 \noindent In this section, we provide the consistency
 of the estimates
(\ref{estimate}). We also state their asymptotic normality and
evaluate their limiting variance.
 Statistics of the form
\[\mathbf{\Psi}_n :=\int_{\mathbf{I}}\psi(\mathbf{u})~d\mathbf{C}_{n}(\mathbf{u}),\] belong to the general
class of \emph{multivariate rank statistics}. Their asymptotic
properties have been investigated at length by a number of
authors, among whom we may cite \cite{Ruymgaart_Shorack_Zwet1972},
\cite{Ruymgaart1974} and \cite{Ruschendorf1976}. In particular,
the previous authors have provided regularity conditions, imposed
on $\psi(\cdot)$, which imply the asymptotic normality of
$\mathbf{\Psi}_n$. The corresponding arguments have been modified by
\cite{Genest_Ghoudi_Rivest1995}, as to establish almost sure
convergence of the estimators that they consider (see, e.g.,
\cite{Genest_Ghoudi_Rivest1995} Proposition A.1). In the same
spirit, the limiting behavior, as $n$ tends to the infinity, of
the estimators and test statistics which we will introduce later
on, will make an instrumental use of the general theory of
multivariate rank statistics, and rely, in particular, on
Proposition A.1 in \cite{Genest_Ghoudi_Rivest1995}. The existence
and consistency of our estimators will be established through an
application of the law of the iterated logarithm for empirical
copula processes, in combination with general arguments from
multivariate rank statistics theory (we refer to
\cite{deheuvels1979a}, \cite{fermanianradulovicdragan2004} and
references therein). We will use the following notations
\begin{equation*}
 \mathds{K}_1(\theta,\mathbf{u}):=\phi^{\prime}
\left(\frac{1}{\mathbf{c}_\theta(\mathbf{u}) }\right)
\end{equation*} and
\begin{equation*}
\mathds{K}_2(\theta,\mathbf{u}):=\left\{\phi^{\prime}\left(\frac{1}{\mathbf{c}_{\theta}(\mathbf{u})}\right)
 \frac{1}{\mathbf{c}_{\theta}(\mathbf{u})}-
 \phi\left(\frac{1}{\mathbf{c}_{\theta}(\mathbf{u})}\right)\right\}.
 \end{equation*}
\begin{definition}
\begin{enumerate}
\item[(i)] Let $\mathscr{Q}$ be the set of continuous functions
$q$ on $[0,1]$ which are positive on $(0,1)$, symmetric about
$1/2$, increasing on $[0,1/2]$ and satisfy
$\int_0^1\{q(t)\}^{-2}dt< \infty$.
 \item[(ii)] A function
$r: (0,1)\rightarrow (0,\infty)$  is called u-shaped if it is
symmetric about $1/2$ and increasing on $(0,1/2]$. \item[(iii)]
For $0<\beta < 1$ and u-shaped function $r$, we define
$$
r_\beta(t)=\left\{
\begin{array}{lcr}
  r(\beta t)& if & 0<t\leq 1/2; \\
   r\{1-\beta(1-t)\}& if & 1/2<t\leq 1/2. \\
\end{array}\right.
$$
If for $\beta > 0$ in a neighborhood of $0$, there exists a
constant $M_\beta$, such that $r_\beta \leq M_\beta r$ on $(0,1)$,
then $r$ is called a reproducing u-shaped function. We denote by
$\mathscr{R}$ the set of reproducing u-shaped functions.
\end{enumerate}
\end{definition}
\noindent Typical examples of elements in $\mathscr{Q}$ and
$\mathscr{R}$ are given by \begin{equation*}
q(t)=\left[t(1-t)\right]^\zeta,~
0<\zeta<1/2,~~~~r(t)=\varrho\left[t(1-t)\right]^{-\varsigma},~\varsigma\geq
0,~\varrho\geq 0.
\end{equation*}
We make use of the following conditions.
\begin{enumerate}
    \item[(C.1)] There exists a neighborhood $N(\theta_\mathbf{T})\subset \mathbf{\Theta}_e$ of $\theta_\mathbf{T}$
     such that the first and the second partial derivatives with respect to
    $\theta$ of $\mathds{K}_1(\theta,\mathbf{u})$ are dominated on $N(\theta_\mathbf{T})$ by some $\lambda$-integrable
functions;
\end{enumerate}
\begin{enumerate}
    \item[(C.2)] There exists a neighborhood $N(\theta_\mathbf{T})$ of $\theta_\mathbf{T}$, such that for all
    $\theta \in N(\theta_\mathbf{T})$, the functions $\frac{\partial}{\partial\theta_i}\,\mathbf{m}(\theta,\mathbf{u}): (0,1)^d\rightarrow
    \mathbf{R}$ are continuously differentiable, and there exist functions $r_i \in
    \mathscr{R}$, $\widetilde{r}_i \in \mathscr{R}$
    and $q \in \mathscr{Q}$ ($i,j=1,\ldots,d,~~i\neq j$ and $\ell,\ell^\prime,\ell^{\prime\prime}=1,\ldots,p$) with
    \begin{enumerate}
    \item[(i)]$$
    \big|\frac{\partial}{\partial\theta_\ell}\,\mathbf{m}(\theta,\mathbf{u})\big|\leq \prod_{i=1}^dr_i(u_i),~~
    \big|\frac{\partial^2}{\partial\theta_\ell\partial u_i}\,\mathbf{m}(\theta,\mathbf{u})\big|\leq
    \widetilde{r}_i(u_i)\prod_{i\neq j}^dr_j(u_j);$$
\item[(ii)] $ \big|\frac{\partial^3}{\partial
\theta_\ell\partial\theta_{\ell^\prime}\partial\theta_{\ell^{\prime\prime}}}\mathds{K}_2(\theta,\mathbf{u})\big|\leq
    \prod_{i=1}^dr_i(u_i);
$ \item[(iii)]$$ \big|\mathbf{m}(\theta,\mathbf{u})\big|\leq
\prod_{i=1}^dr(u_i),  \big|\frac{\partial}{\partial
u_i}\mathbf{m}(\theta,\mathbf{u})\big|\leq
\widetilde{r}_i(u_i)\prod_{i\neq j}^dr_j(u_j);
   $$
\item[(iv)]$$\big|\frac{\partial
}{\partial\theta_\ell}\mathbf{m}(\theta,\mathbf{u})\big|^2\leq\prod_{i=1}^dr_i(u_i),~~
\big|\frac{\partial^2
}{\partial\theta_\ell\partial\theta_{\ell^\prime}}\mathbf{m}(\theta,\mathbf{u})\big|\leq
    \prod_{i=1}^dr_i(u_i)$$
 \end{enumerate}
    and
$$ \int_{\mathbf{I}}\left\{\prod_{i=1}^dr_i(u_i)\right\}^2
d\mathbf{C}_{\theta_\mathbf{T}}(\mathbf{u})< \infty, $$ $$
\int_{\mathbf{I}}\left\{q_i(u_i)\widetilde{r}_i(u_i)\prod_{i\neq
j}^dr_j(u_j)\right\} d\mathbf{C}_{\theta_\mathbf{T}}(\mathbf{u})<
\infty, ~~\mbox{for}~~ i=1,\ldots,d;$$ \item[(C.3)] The matrix
$\int_{\mathbf{I}} (\partial^2/
\partial^2\theta)
\mathbf{m}(\theta,\mathbf{u})d\mathbf{C}_{\theta_\mathbf{T}}(\mathbf{u})$
is non singular;

\item[(C.4)] The function $\mathbf{u}\in I\mapsto
\frac{\partial}{\partial\theta}\,
\mathbf{m}(\theta_\mathbf{T},\mathbf{u})$ is of bounded variation
on $I$.
\end{enumerate}
The main result to be proved here may now be stated precisely as
follows.
\begin{theorem}\label{theoreme1}
Assume that conditions C.1-C.4 hold.
\begin{enumerate}
\item Let $B(\theta_{\mathbf{T}},n^{-1/3}):=\left\{\theta \in
\mathbf{\Theta}_{e} , \|\theta -\theta_{\mathbf{T}} \| \leq
n^{-1/3} \right\}$, then as $n$ tends to infinity, with
probability one, the function $\theta\mapsto\int_{\mathbf{I}}
\mathbf{m}(\theta,\mathbf{u})~d\mathbf{C}_{n}(\mathbf{u})$ attains
its
 maximum value at some point $\widehat\theta_n$ in the interior of  $B(\theta_{\mathbf{T}},n^{-1/3})$, which implies that the estimate
 $\widehat\theta_n$ is consistent and satisfies
  \begin{equation*}
  \int_{\mathbf{I}} \frac{\partial}{\partial \theta}\mathbf{m}(\widehat\theta_n,\mathbf{u})~d\mathbf{C}_{n}(\mathbf{u}) = 0.
 \end{equation*}
\item  $\sqrt{n}(\widehat\theta_n -\theta)$ converges in
distribution to a centered multivariate normal random variable
with covariance matrix
\begin{equation}\label{fin}
\mathbf{\Xi}_\phi=\mathbf{S}^{-1}\mathbf{M}\mathbf{S}^{-1},
\end{equation}
with
\begin{equation}\label{S}
\mathbf{S}:=-\int_{\mathbf{I}} \frac{\partial^2} {\partial
\theta^2}
\mathbf{m}(\theta_\mathbf{T},\mathbf{u})d\mathbf{C}_{\theta_\mathbf{T}}(\mathbf{u}),
\end{equation}
 and
\begin{equation}\label{M}
\mathbf{M}:=\mathbf{Var}\left[\frac{\partial}{\partial
\theta}\mathbf{m}(\theta_\mathbf{T},F_1(X_1),\ldots,F_d(X_d))+\sum_{i=1}^d\mathds{W}_{i}(\theta_\mathbf{T},X_i)\right],
\end{equation}
where
\begin{equation*}
 \mathds{W}_{i}(\theta_\mathbf{T},X_i) := \int_{\mathbf{I}} \left\{\mathds{1}_{\{F_{i}(X_i) \leq
 u_i\}}-u_i\right\}{\textstyle\frac{\partial^2}{\partial\theta \partial
 u_i}}\mathbf{m}\left(\theta_\mathbf{T},\mathbf{u}\right)
 d\mathbf{C}_{\theta_\mathbf{T}}(\mathbf{u}),~i=1,\ldots,d.
\end{equation*}
\end{enumerate}
\end{theorem}
The proof of Theorem \ref{theoreme1} is postponed to the Appendix.
\begin{remark}
The aim of Theorem \ref{theoreme1} part (a) is not to establish
the optimal rate of the estimate  but merely  the existence and
the consistency (a.s.) of the estimate. We have considered
$n^{1/3}$ because it works well, indeed, in Taylor expansion
(\ref{dev 1}), in the proof, the third term of the RHS is $O(1)$
only for this rate, which is the major key of the demonstration.
\end{remark}
\section{New tests of independence}
 \noindent One of our
motivation is to build a statistical test of independence, based
on $\phi$-divergence. In the framework of the parametric copula
model, the null hypothesis, i.e., the independence case
$$\mathbf{C}_{\theta_0}(u_1,\ldots,u_d) = \prod_{i=1}^du_i$$
corresponds to $$\mathscr{H}_0~:~\theta_\mathbf{T} = \theta_0.$$  We
consider the composite alternative  hypothesis
$$\mathscr{H}_1~:~\theta_\mathbf{T} \neq \theta_0.$$ Since,
$\theta_0$ is a boundary value of the parameter space
$\mathbf{\Theta}$, we can see that the convergence in distribution
of the corresponding pseudo-likelihood ratio statistic to a
$\chi^2$ random variable does not hold; see
\cite{bouzebda-keziou2008}. We give now a solution to this
problem. We propose the following statistics
\begin{equation}
\mathbf{T}_n:=\frac{2n}{\phi^{\prime
\prime}(1)}\widehat{D_\phi}(\theta_0,\theta_\mathbf{T}).
\end{equation}
We will see that the proposed statistic converges in distribution,
under the null hypothesis $\mathscr{H}_0$, to a $\chi^2$ random
variable with $p$ degrees of freedom, which permits to build a
test of $\mathscr{H}_0$ against $\mathscr{H}_1$ asymptotically of
level $\alpha$. The limit law of $\mathbf{T}_n$ is given also
under the alternative hypothesis $\mathscr{H}_1$. We will use the
following additional conditions.
\begin{enumerate}
\item [(C.5)] We have
$$\lim_{\theta\rightarrow \theta_0 }\frac{\partial^2}{\partial\theta_\ell
\partial u_i}\,\mathbf{m}(\theta,\mathbf{u})=0,$$ and there exists a neighborhood $N(\theta_0)$ of $\theta_0$ and there exist functions $r_i \in
    \mathscr{R}$, $\widetilde{r}_i \in \mathscr{R}$
    and $q_i \in \mathscr{Q}$ ($i=1,\ldots,d$ and $\ell=1,\ldots,p$), such that for all
    $\theta \in N(\theta_0)$,
\begin{equation*}
\big|\frac{\partial^2}{\partial\theta_\ell\partial
u_i}\,\mathbf{m}(\theta,\mathbf{u})\big|<\widetilde{r}_i(u_i)\prod_{j\neq
j}^dr(u_j)
\end{equation*}
and
$$
\int_{\mathbf{I}}\left\{q_i(u_i)\widetilde{r}_i(u_i)\prod_{i\neq
j}^dr_j(u_j)\right\} d\mathbf{C}_{\theta_\mathbf{T}}(\mathbf{u})<
\infty.$$
 \end{enumerate}
 \begin{remark}
When $\theta_\mathbf{T}=\theta_0$, under the conditions (C.1) and
(C.5) we can see that $\mathbf{S}$ and $\mathbf{M}$ can be written
as
$$\mathbf{S}=\mathbf{M}=\int_{\mathbf{I}} \left[\frac{\partial} {\partial \theta}
\mathbf{m}(\theta_\mathbf{T},\mathbf{u})\right]\left[\frac{\partial}
{\partial \theta}
\mathbf{m}(\theta_\mathbf{T},\mathbf{u})\right]^\top
d\mathbf{C}_{\theta_\mathbf{T}}(\mathbf{u}).$$
\end{remark}
The following theorem gives the limiting law of the statistics
$\mathbf{T}_n$ under the both hypothesis $\mathscr{H}_0$ and
$\mathscr{H}_1$.
\begin{theorem}\label{theoreme2}
\begin{enumerate}
\item [(1)] Assume that conditions C.1-C.5 hold. If
$\theta_\mathbf{T}=\theta_{0}$, then the statistic $\mathbf{T}_n$
converges in distribution to a $\chi^{2}$ variable with $p$
degrees of freedom. \item [(2)] Assume that conditions C.1-C.5
hold. If $\theta_\mathbf{T}\neq \theta_{0}$, then
$$\sqrt{n}\left(\widehat{D_\phi}(\theta_0,\theta_{\mathbf{T}})-D_\phi(\theta_0
,\theta_{\mathbf{T}})\right)$$ converges in distribution to a
centered normal variable with variance
\begin{equation}\label{mm}
\sigma^2_{\phi}(\theta_0,\theta_\mathbf{T}) :=
\mathbf{Var}\left[\mathbf{m}(\theta_\mathbf{T},F_1(X_1),\ldots,F_d(X_d))+\sum_{i=1}^d\mathds{Y}_{i}(\theta_\mathbf{T},X_i)
\right],
\end{equation}
\end{enumerate}
where
\begin{equation*}
 \mathds{Y}_{i}(\theta_\mathbf{T},X_i) := \int_{\mathbf{I}} \left\{\mathds{1}_{\{F_{i}(X_i) \leq
 u_i\}}-u_i\right\}{\textstyle\frac{\partial}{ \partial
 u_i}}\mathbf{m}\left(\theta_\mathbf{T},\mathbf{u}\right)
 \mathbf{c}_{\theta_\mathbf{T}}(\mathbf{u})~du_1\ldots du_d,~i=1,\ldots,d.
\end{equation*}
\end{theorem}
The proof of Theorem \ref{theoreme2} is postponed to the Appendix.
\begin{remark} An application of Theorem \ref{theoreme2}, leads to reject the
null hypothesis $\mathscr{H}_0:\theta_\mathbf{T}=\theta_0$,
whenever the value of the statistic $ \mathbf{T}_n$ exceeds
$q_{1-\alpha}$, namely, the $(1-\alpha)$-quantile of the $\chi^2$
law with $p$ degrees of freedom. The corresponding test is then,
asymptotically of level $\alpha$, when $n\rightarrow\infty$. The
critical region is, accordingly, given by
\begin{equation*}
CR := \left\{{\bf T}_n > q_{1-\alpha}\right\}.
\end{equation*}
The fact that this test is consistent follows from Theorem
\ref{theoreme2}. Further, this theorem can be used to give an
approximation to the power function
$\theta_\mathbf{T}\in\mathbf{\Theta}\mapsto
\beta(\theta_\mathbf{T}):=P_{\theta_\mathbf{T}}\left\{CR\right\}$
in a similar way to  \cite{KeziouLeoni2007}. We so obtain that
\begin{equation}\label{marre2}
\beta(\theta_\mathbf{T})\approx 1-\Phi\left(
\frac{\sqrt{n}}{\sigma_{\phi}(\theta_0,\theta_\mathbf{T})}
\left(\frac{q_{1-\alpha}}{2n}-D_\phi(\theta_0
,\theta_{\mathbf{T}})\right)\right),
\end{equation}
where $\Phi$ denotes, as usual, the cumulative distribution
function of a $\mathcal{N}(0,1)$ standard normal random variable. A useful
consequence of (\ref{marre2}) is the possibility of computing an
approximate value of the sample size ensuring a specified power
$\beta(\theta_\mathbf{T})$, with respect to some pre-assigned
alternative $\theta_\mathbf{T}\neq\theta_0.$ Let $n_0$ be the
positive root of the equation
\begin{equation*}
\beta = 1-\Phi\left(
\frac{\sqrt{n}}{\sigma_{\phi}(\theta_0,\theta_\mathbf{T})}
\left(\frac{q_{1-\alpha}}{2n}-D_\phi(\theta_0
,\theta_{\mathbf{T}})\right)\right),
\end{equation*}
which can be rewritten into
\begin{equation*}
n_0=\frac{(a+b)- \sqrt{a(a+2b)}}{2D_\phi(\theta_0
,\theta_{\mathbf{T}})^2},
\end{equation*}
where $a
:=\sigma_{\phi}(\theta_0,\theta_\mathbf{T})\left(\Phi^{-1}(1-\beta)\right)^2$
and $b:=q_{1-\alpha} D_\phi(\theta_0 ,\theta_{\mathbf{T}})$. The
sought-after approximate value of the sample size is then given
$$n^*:=\lfloor n_0\rfloor+1,$$
where $\lfloor u\rfloor$ denote the integer part of $u$.
\end{remark}
\begin{remark}
For point estimation, the estimator based
on $\phi$-divergence when we extend the parameter space, may not have a meaningful interpretation and most probably has a larger mean square error.
However, from Theorem \ref{theoreme1} and \ref{theoreme2}, it is clear that
an asymptotic $1-\alpha$ confidence interval or region,
$\mathscr{R}_\alpha$ about $\theta$ can be easily constructed
using the intersection method as described in \cite{feng1992}.
\end{remark}

\begin{remark}
The above regularity conditions are satisfied for a large number of single-parameter
families of bivariate copulas including the standard bivariate normal, the Farlie-Gumbel-Morgenstern
system, and copulas of the Archimedean variety such as those of Ali-Mikhail-Haq and Frank;
see, e.g., \cite{Genest_Ghoudi_Rivest1995} and \cite{tsukahara2005}.  Note that the score functions for some copulas are
unbounded near the origin or the point $(1,1)$, so we need to know the above regularity
 conditions, at least theoretically as be mentioned in \cite{tsukahara2005}.
\end{remark}

\begin{remark}
The parameters   (\ref{S}) and (\ref{M})   may be consistently
estimated respectively  by the sample mean of
\begin{equation}\label{S1}
 \frac{\partial^2}{\partial \theta^2}
\mathbf{m}(\widehat{\theta}_n,F_{1n}(X_{1,k}),\ldots,F_{dn}(X_{d,k})),
\quad k=1,\ldots,n,
\end{equation}
and the sample variance of
\begin{equation}
 {\textstyle\frac{\partial}{\partial\theta}}\,\mathbf{m}\left(\widehat{\theta}_n,F_{1n}
 (X_{1,k}),\ldots,F_{dn}(X_{d,k})\right) + \sum_{i=1}^d\mathds{W}_{i}(\widehat{\theta}_n,X_{i,k}), \quad k=1,\ldots,n,
\end{equation}
as was done in \cite{Genest_Ghoudi_Rivest1995}. The asymptotic
variance (\ref{mm}) can be consistently estimated in the same way.
\end{remark}
\begin{remark} The set $\mathbf{\Theta}_e$ defined in (\ref{def de Theta_e }) is generally with non empty interior $\mathring{\mathbf{\Theta}}_e$.
In particular, we may check that $\theta_0$ (the value
corresponding to independence) belongs to
$\mathring{\mathbf{\Theta}}_e$, since the integral in (\ref{def de
Theta_e }) is finite; it is equal to zero when $\theta =
\theta_0$, for any copula density $\mathbf{c}_\theta(\cdot)$.
However, it is hard to determine the whole set $\mathbf{\Theta}_e$
for some copulas, but in order to test the independence, we need
only to prove the existence of a neighborhood $N(\theta_0)$ of
$\theta_0$ for which the integral in (\ref{def de Theta_e }) is
finite since we calculate the estimate $\widehat{\theta}_n$ in
(\ref{estimate}) by Newton-Raphson algorithm using $\theta_0$ as
initial point. The explicit calculation of the integral in
(\ref{def de Theta_e }) may be complicated for some copulas, in
such cases we use the Monte Carlo method to compute this integral.
\end{remark}
\section{Simulations}\label{simula}
In this section, we report the results from simulation experiments
carried out to assess the performance of the proposed estimators.
To this end, we have considered the FGM copula. For the experiment
considered here, we compute the MPL, KL-divergence, $\chi^2$-divergence, Hellinger divergence and some power divergence estimates,
 and report their variance value, bias  and mean-squared error.
In order to compare the robustness of the proposed estimates we
consider several scenarios of contamination. To be more precise,
we considered $\epsilon$-contaminated models, where a proportion
$\epsilon$ of observations were replaced by atypical ones
generated from a contaminating distribution
$\mathbf{F}^*(\cdot,\cdot)$. We set $\epsilon$ equal to $0\%$,
$5\%$, and $10\%$, and $\mathbf{F}^*(\cdot,\cdot)$  as the
bivariate normal distribution with correlation coefficient $\rho=
0.00$ and very small variances, acting as a point mass
contamination as in \cite{Mendes2007}. The sample size is $n=500$
and the estimates are obtained from $1000$ independent runs. \\

\noindent (i) Under no contamination: All procedures showed
reasonable accuracy. The Hellinger and $\chi^2$-divergence
estimates seem to be as good as the MPL estimator. This is more
evident as
the sample size gets larger; see  Table \ref{tabnoncon}. \\

\noindent (ii) Under $5\%$ and $10\%$ contamination: MPL estimator
is recommended when there is no contamination but its
performance may deteriorate rapidly if the sample is pooled, see
Tables \ref{tabcon1} and \ref{tabcon2}. In the contaminated case 
the power divergence estimator with $\gamma=2.5$ is superior with
respect to the others. It seems that the $\chi^2$-divergence
estimate behaves well also
for contaminated data. 
\\

\noindent From the three  tables  \ref{tabnoncon}, \ref{tabcon1}
and \ref{tabcon2}, we can see that the  $\chi^2$-divergence
estimates is a good trad-off between efficiency and robustness.\\

\noindent  In future work, it would be interesting to provide a
complete investigation  of robustness of semiparametric copula
estimator which requires nontrivial mathematics, this would go
well beyond the scope of the present paper.

\begin{table}[h!]
\centering
\begin{tabular}{|l|c|c|c|c|}
\hline
\multicolumn{5}{|c|}{$\theta_T=0.01$,~~$n=500$,~~$rep=1000$ }  \\
\hline Divergence
& Estimate& Variance & Bias&MSE\\
\hline $\gamma = - 0.5$ &

    0.0924  &  0.0177  &  0.0076 &   0.0177\\
\hline

$\gamma =0$ (MPL) &

    0.0920   & 0.0176 &   0.0080  &  0.0176\\
\hline
$\gamma =0.5$ (H)&

    0.0916 &   0.0175 &   0.0084   & 0.0176\\
\hline

 $\gamma=1$ (KL)&

    0.0913  &  0.0175 &   0.0087  &  0.0176\\
\hline

$\gamma = 1.5$ &

    0.0911  &  0.0176  &  0.0089   & 0.0177\\
\hline

$\gamma =2$ ($\chi^2$)  &

    0.0911 &   0.0178 &   0.0089   & 0.0179\\
\hline

$\gamma=2.5$ &

    0.0912  &  0.0181  &  0.0088  &  0.0181\\

    \hline
\end{tabular}
\caption{No contamination: $\epsilon= 0.00$.~~}
\label{tabnoncon}
\end{table}

\begin{table}[h!]
\centering
\begin{tabular}{|l|c|c|c|c|}
\hline
\multicolumn{5}{|c|}{$\theta_T=0.01$,~~$n=500$,~~$rep=1000$ }  \\
\hline Divergence
& Estimate& Variance & Bias&MSE\\
\hline
  $\gamma =-0.5$ &

    0.0949  &  0.0191 &   0.0051  &  0.0191\\
\hline

$  \gamma =0$ (MPL) &

     0.0921  &  0.0181 &   0.0079 &   0.0181\\
\hline

  $\gamma=0.5$ (H)&

     0.0895   & 0.0172  &  0.0105  &  0.0173\\
\hline

 $\gamma =1$ (KL)&

      0.0873  &  0.0164  &  0.0127   & 0.0166\\
\hline

 $\gamma =1.5$ &

     0.0853   & 0.0158   & 0.0147&    0.0160\\
\hline

  $\gamma=2$ ($\chi^2 $)&

    0.0835  &  0.0153 &   0.0165 &   0.0155\\
\hline

  $\gamma=2.5$ &

   0.0820   & 0.0149 &   0.0180  &  0.0152\\
    \hline
\end{tabular}
\caption{Contamination: $\epsilon= 0.05$.~~}
\label{tabcon1}
\end{table}
    
\begin{table}[h!]
\centering
\begin{tabular}{|l|c|c|c|c|}
\hline
\multicolumn{5}{|c|}{$\theta_T=0.01$,~~$n=500$,~~$rep=1000$ }  \\
\hline Divergence
& Estimate& Variance & Bias&MSE\\
\hline
    $\gamma=-0.5$ &

    0.0916  &  0.0191   & 0.0084  &  0.0191\\

  \hline
  $\gamma =0$ (MPL)&

    0.0867  &  0.0171 &   0.0133 &   0.0173\\

  \hline
 $\gamma=0.5$ (H)&

     0.0825   & 0.0155  &  0.0175  &  0.0158\\

  \hline
 $\gamma =1$ (KL)&

    0.0787    &0.0141 &   0.0213 &   0.0146\\

  \hline
  $\gamma=1.5$ &

    0.0754    &0.0130   & 0.0246   & 0.0136\\

  \hline
  $ \gamma=2$ ($\chi^2$) &

    0.0724    &0.0120  &  0.0276  &  0.0128\\

  \hline
 $\gamma= 2.5$ &

    0.0698    &0.0112    &0.0302    &0.0121\\
    \hline
\end{tabular}
\caption{Contamination: $\epsilon= 0.10$.~~}
\label{tabcon2}
\end{table}
\section{Concluding remarks}
\noindent We have introduced a new estimation and test procedure
in parametric copula models with unknown margins. The method is
based on divergences between copulas and the duality technique. It
generalizes the maximum pseudo-likelihood one, and applies both
when the parameter is an interior or a boundary value, in
particular for testing the null hypothesis of independence.
Simulation results show that the $\chi^2$-divergence estimate is a
good trad-off between efficiency and robustness.  It will be
interesting to investigate theoretically the problem of the choice
of the divergence which leads to an ``\emph{optimal}'' (in some sense)
estimate or test in terms of efficiency and robustness.

\section{Appendix}\label{appendix} \noindent First we give a technical Lemma which
we will use to prove our results.
\begin{lemma}\label{lemma1}
Let $\mathbf{F}_{\theta_\mathbf{T},F_1,\ldots,F_d}(\cdot)$ have a
continuous margins and let $\mathbf{C}_{\theta_\mathbf{T}}(\cdot)$
have continuous partial derivatives. Assume that
$\mathbf{\mathbf{\xi}}(\cdot)$ is a continuous function, with
bounded variation. Then
\begin{equation}
\int_{\mathbf{I}} \mathbf{\xi}(\mathbf{u})~d\left(\mathbf{C}_{n}(\mathbf{u})- \mathbf{C}(\mathbf{u})\right)=
O\left(n^{-1/2}(\log\log n)^{1/2}\right)~~ (a.s.).
\end{equation}
\end{lemma}

\noindent{\bf{Proof of Lemma \ref{lemma1}.}}  Recall that the
\emph{modified empirical copula} $\mathbf{C}_n(\cdot)$, is slightly different from
the \emph{empirical copula} $\mathds{C}_n(\cdot)$, introduced by
\cite{deheuvels1979a}, and defined by
\begin{equation}\label{marre3}\mathds{C}_n(\mathbf{u})=\mathbf{F}_n\Big(F_{1n}^{-1}(u_1),
\ldots,F_{dn}^{-1}(u_d)\Big)\quad\hbox{for}\quad\mathbf{u}\in
(0,1)^d,\end{equation} where $F_{in}^{-1}(\cdot)$ for
$i=1,\ldots,d$ denote the empirical quantile functions, associated
with $F_{in}(\cdot)$ for $i=1,\ldots,d$, respectively, and defined
by
\begin{equation*}
F_{in}^{-1}(t):=\inf\{x\in \mathbf{R}\mid F_{in}(x)\geq t\},\quad
i=1,\ldots,d.
\end{equation*}
Note that the subtle difference lies in the fact that $\mathds{C}_n(\cdot)$ is
left-continuous with right-hand limits, whereas $\mathbf{C}_n(\cdot)$ on the other
hand is right continuous with left-hand limits. The difference
between $\mathds{C}_n(\cdot)$ and $\mathbf{C}_n(\cdot)$, however, is small
\begin{equation}\label{app 1}
\sup_{{\bf{u}}\in
 I}\left|{\mathds{C}}_{n}({\bf{u}})-\mathbf{C}_{n}({\bf{u}})\right|=\frac{1}{n}.
\end{equation}
As in the proof of Lemma 5.1 in
\cite{bouzebda-keziou2008}, using integration by parts, we can
prove that there exists a constant $\kappa>0$, depending upon $d$
only, such that
\begin{equation*}
\left|\sqrt{n}\int_{\mathbf{I}}\mathbf{\xi}(\mathbf{u})~d(\mathbf{C}_{n}-\mathbf{C})(\mathbf{u})\right|\leq
\kappa \sqrt{n}\sup_{{\bf u}\in I}\left|(\mathbf{C}_{n} - \mathbf{C})(\bf
u)\right|\int_{\mathbf{I}}d\left|\mathbf{\xi}(\bf u)\right|.
\end{equation*}
Or, by Fubini's Theorem, we can write
\begin{eqnarray*}
\lefteqn{\left|\sqrt{n}\int_{\mathbf{I}}\mathbf{\xi}(u_1,\ldots,u_d)~d(\mathbf{C}_{n}-\mathbf{C})(u_1,\ldots,u_d)\right|}\\
&=&\left|\sqrt{n}\int_{\mathbf{I}}\left\{\int_0^{u_1}\cdots\int_0^{u_d}d\mathbf{\xi}(s_1,\ldots,s_d)
\right\}~d(\mathbf{C}_{n}-\mathbf{C})(u_1,\ldots,u_d)\right|
\\&=&\left|\sqrt{n}\int_{\mathbf{I}}\left\{\int_\mathbf{I}\mathds{1}_{\{s_1\leq u_1\}}\ldots
\mathds{1}_{\{s_d\leq
u_d\}}d\mathbf{\xi}(s_1,\ldots,s_d)\right\}~d(\mathbf{C}_{n}-\mathbf{C})(u_1,\ldots,u_d)\right|
\\&=&\left|\sqrt{n}\int_{\mathbf{I}}\left\{\int_\mathbf{I}\mathds{1}_{\{s_1\leq u_1\}}\ldots
\mathds{1}_{\{s_d\leq
u_d\}}d(\mathbf{C}_{n}-\mathbf{C})(u_1,\ldots,u_d)\right\}~d\mathbf{\xi}(s_1,\ldots,s_d)\right|
\\&=&\left|\sqrt{n}\int_{\mathbf{I}}\left\{\Delta_\mathbf{s}^\mathbf{1}(\mathbf{C}_{n}-
\mathbf{C})(u_1,\ldots,u_d)\right\}~d\mathbf{\xi}(s_1,\ldots,s_d)\right|
\\&\leq&(2^d-1)\sqrt{n}\sup_{{\bf u}\in I}\left|(\mathbf{C}_{n} - \mathbf{C})(\bf
u)\right|\int_{\mathbf{I}}d\left|\mathbf{\xi}(\bf u)\right|,
\end{eqnarray*}
where for $\mathbf{a}$ and $\mathbf{b}$ in $\mathbf{I}$
\begin{equation*}
\Delta_\mathbf{a}^\mathbf{b}(\mathbf{C}_{n}-\mathbf{C})(\mathbf{u}):=\Delta_{a_d}^{b_d}\Delta_{a_{d-1}}^{b_{d-1}}\cdots
\Delta_{a_2}^{b_2}\Delta_{a_1}^{b_1}(\mathbf{C}_{n}-\mathbf{C})(\mathbf{u})
\end{equation*}
and for $j=1,\ldots,d$
\begin{eqnarray*}
\Delta_{a_j}^{b_j}(\mathbf{C}_{n}-\mathbf{C})(\mathbf{u})&:=&(\mathbf{C}_{n}-\mathbf{C})(u_1,\ldots,u_{j-1},b_j,u_{j-1},\ldots,u_d)\\
&&-(\mathbf{C}_{n}-\mathbf{C})(u_1,\ldots,u_{j-1},a_j,u_{j-1},\ldots,u_d).
\end{eqnarray*}
One may check (see  Theorem 3.1 in \cite{deheuvels1979a}) that
there exists a constant $\gamma$ (depending upon
$\mathbf{C}(\cdot)$ only) such that, with probability $1$,
\begin{equation}\label{deheuv}
\limsup_{n\rightarrow \infty}\left\{\frac{n}{\log \log n}\right\}^{1/2}\sup_{\mathbf{u}\in I}
\left|\mathds{C}_n(\mathbf{u})-\mathbf{C}(\mathbf{u})\right|=\gamma<\infty.
\end{equation}
From this and (\ref{app 1}), applying (\ref{deheuv}), we obtain
$$\int_{\mathbf{I}} \mathbf{\xi}(\mathbf{u})~d(\mathbf{C}_{n}-\mathbf{C})(\mathbf{u})=
O\left(n^{-1/2}(\log\log n)^{1/2}\right)~~ (a.s.).$$ \hfill $\blacksquare$\\
\noindent\textbf{Proof of Theorem \ref{theoreme1}} (1) Under the
Assumptions (C.1) and (C.2.ii), a simple calculation gives
\begin{equation}\label{5.3}
\int_{\mathbf{I}}\frac{\partial}{\partial \theta}\mathbf{m}(\theta
,\mathbf{u}) d\mathbf{C}_{\theta_\mathbf{T}}(\mathbf{u}) = 0,
\end{equation}
and
\begin{equation}\label{vv}
\int_{\mathbf{I}}\frac{\partial^2}{\partial
\theta^2}\mathbf{m}(\theta ,\mathbf{u})
d\mathbf{C}_{\theta_\mathbf{T}}(\mathbf{u})=-\int_{\mathbf{I}}\phi
^{\prime\prime}\left(\frac{1}{\mathbf{c}_{\theta_\mathbf{T}}}\right)
\frac{{\dot{\mathbf{c}}_{\theta_\mathbf{T}}
\dot{\mathbf{c}}_{\theta_\mathbf{T}
}^{\top}}}{\mathbf{c}_{\theta_\mathbf{T}}^3}~d\lambda=-\mathbf{S}.
\end{equation}
We see that the matrix $\mathbf{S}$ is symmetric and positive
using the fact that the second derivative $\phi^{''}(\cdot)$ is
nonnegative by the assumption that the function $\phi(\cdot)$ is
convex. Hence, $\mathbf{S}$ is positive definite by (C.3).
Introduce the statistic $\mathbf{\Phi}_n(\theta_\mathbf{T})$
defined by
\begin{equation}
\mathbf{\Phi}_n(\theta_\mathbf{T}):=   \int_{\mathbf{I}}
\frac{\partial}{\partial
\theta}\mathbf{m}(\theta,\mathbf{u})~d\mathbf{C}_{n}(\mathbf{u}),
\end{equation}
and combine (\ref{5.3}) and condition (C2)(i) with Theorem 2.1 in
\cite{Ruymgaart_Shorack_Zwet1972} to show that, as $n\rightarrow
\infty $
\begin{equation}\label{nom}
\sqrt{n}\mathbf{\Phi}_n(\theta_\mathbf{T})\stackrel{d}{\rightarrow}
\mathcal{N}(0,\mathbf{M}),
\end{equation}
where $\mathbf{M}$ is defined in (\ref{M}).
We can refer also to the
Proposition 3 in \cite{tsukahara2005} for the same result in
(\ref{nom}).
 Denote
\begin{equation}\label{vas}
\mathbf{\Upsilon}_n(\theta_\mathbf{T}):=\int_{\mathbf{I}}
\frac{\partial^2}{\partial
\theta^2}\mathbf{m}(\theta,\mathbf{u})~d\mathbf{C}_{n}(\mathbf{u}),
\end{equation}
we make use of (\ref{vv}) and (C.2.iv) in connection with
Proposition A.1 of \cite{Genest_Ghoudi_Rivest1995}, one finds
\begin{equation}\label{saa}
\mathbf{\Upsilon}_n(\theta_\mathbf{T})\rightarrow -\mathbf{S},~~
(a.s.).
\end{equation}
We recall that $\mathbf{S}$ is in (\ref{S}). Now, for any $\theta
= \theta_\mathbf{T} + \textbf{v} n^{-1/3}$ with $\|\textbf{v}\|
\leq 1$, consider a Taylor expansion of  $\int_{\mathbf{I}}
\mathbf{m}(\theta,\mathbf{u})~d\mathbf{C}_{n}(\mathbf{u})$ in
$\theta$ around $\theta_\mathbf{T}$, and use (\ref{saa}), and
(C.2.ii) to obtain
\begin{eqnarray}\label{dev 1}
\nonumber &&n\int_{\mathbf{I}}
\mathbf{m}(\theta,\mathbf{u})d\mathbf{C}_{n}(\mathbf{u})
-n\int_{\mathbf{I}} \mathbf{m}(\theta_{\mathbf{T}},\mathbf{u})~d\mathbf{C}_{n}(\mathbf{u})=\\
 & &
n^{2/3}\textbf{v}^\top\mathbf{\Phi}_n(\theta_\mathbf{T})+2^{-1}n^{1/3}\textbf{v}\mathbf{S}\textbf{v}^\top+O(1)~~(a.s.)
\end{eqnarray}
uniformly in $\textbf{v}$ with $\|\textbf{v}\|\leq 1$. On the
other hand, since $$\int_{\mathbf{I}}
\frac{\partial}{\partial\theta_\ell}\mathbf{m}(\theta_\mathbf{T},\mathbf{u})^{2}~d\mathbf{C}_{\theta_\mathbf{T}}(\mathbf{u})<
\infty,$$ and
$\frac{\partial}{\partial\theta_\ell}\mathbf{m}(\theta,\cdot)$ is
of bounded variation by assumption (C.4)($\ell=1,\ldots,p$), using
Lemma \ref{lemma1} we can show that
\begin{equation}\label{22}
\int_{\mathbf{I}}
\frac{\partial}{\partial\theta_\ell}\mathbf{m}(\theta_\mathbf{T},\mathbf{u})~d\mathbf{C}_{n}(\mathbf{u})=
O\left(n^{-1/2}(\log\log n)^{1/2}\right)~~ (a.s.).
\end{equation}
  Therefore,
using (\ref{dev 1}) and (\ref{22}), we obtain for any $\theta =
\theta_\mathbf{T} + \textbf{v} n^{-1/3}$  with $\|\textbf{v}\|=1$:
\begin{eqnarray}\label{vie}
\nonumber&&n\int_{\mathbf{I}}
\mathbf{m}(\theta,\mathbf{u})~d\mathbf{C}_{n}(\mathbf{u})
-n\int_{\mathbf{I}}
\mathbf{m}(\theta_{\mathbf{T}},\mathbf{u})~d\mathbf{C}_{n}(\mathbf{u})
\\&& \leq O(n^{1/6}(\log\log n)^{1/2})- 2^{-1}\vartheta n^{1/3}
+O(1) ~~(a.s.),
\end{eqnarray}
where $\vartheta$ is the smallest eigenvalue of the matrix
$\mathbf{S}$. Observe  that $\vartheta$ is positive since
$\mathbf{S}$ is symmetric, positive and non singular by assumption
(C.3). Using (\ref{vie}) and the fact that the function $\theta
\mapsto \int_{\mathbf{I}}
\mathbf{m}(\theta_\mathbf{T},\mathbf{u})~d\mathbf{C}_{n}(\mathbf{u})$
is continuous, we conclude that as $n \rightarrow \infty$, with
probability one, $\theta \mapsto \int_{\mathbf{I}}
\mathbf{m}(\theta_\mathbf{T},\mathbf{u})~d\mathbf{C}_{n}(\mathbf{u})$
reaches its maximum value at some point $\widehat \theta_n$
fulfills $$\int_{\mathbf{I}}
\frac{\partial}{\partial\theta}\mathbf{m}(\widehat\theta_n,\mathbf{u})~d\mathbf{C}_{n}(\mathbf{u})=0$$
 and $$\|\widehat\theta_n-\theta_\mathbf{T}\|=O(n^{-1/3}).$$

\noindent(2) Using the first part of Theorem \ref{theoreme1}, by a
Taylor expansion of
$$\int_{\mathbf{I}}
\frac{\partial}{\partial\theta}\mathbf{m}(\widehat\theta_n,\mathbf{u})~d\mathbf{C}_{n}(\mathbf{u}),$$
in $\widehat\theta_n$ around $\theta_\mathbf{T}$, we obtain
\begin{eqnarray*}
\lefteqn{0 = \int_{\mathbf{I}} \frac{\partial}{\partial\theta}
\mathbf{m}(\widehat\theta_n,\mathbf{u})~d\mathbf{C}_{n}(\mathbf{u})}\\
& = & \int_{\mathbf{I}}
\frac{\partial}{\partial\theta}\mathbf{m}(\theta_\mathbf{T},\mathbf{u})~d\mathbf{C}_{n}(\mathbf{u})
 + (\widehat\theta_n - \theta_{\mathbf{T}})^\top\int_{\mathbf{I}}
\frac{\partial^2}{\partial\theta^2}\mathbf{m}(\theta_\mathbf{T},\mathbf{u})~d\mathbf{C}_{n}(\mathbf{u})+o(n^{-1/2}).
\end{eqnarray*}
Hence,
\begin{equation}\label{equivaassmpto}
\sqrt{n}(\widehat\theta_n - \theta_{\mathbf{T}}) =-\left[
\mathbf{\Upsilon}_n(\theta_\mathbf{T})\right]^{-1}\sqrt{n}\mathbf{\Phi}_n(\theta_\mathbf{T})
+o_P(1).
\end{equation}
Using (\ref{nom}) and (\ref{saa}), by Slutsky theorem, we conclude
then
\begin{equation}
\sqrt{n}(\widehat\theta_n - \theta_{\mathbf{T}}) \rightarrow
\mathcal{N}(0,\mathbf{\Xi}_\phi),
\end{equation}
where we recall that $\mathbf{\Xi}_\phi$ is defined in
(\ref{fin}).\hfill$\blacksquare$  \vskip5pt

\noindent\textbf{Proof of Theorem \ref{theoreme2}}\\\noindent(1)
Assume that $\theta_\mathbf{T}=\theta_0$. Hence, from
(\ref{equivaassmpto}), using (\ref{vv}), we obtain
\begin{equation}\label{eqn 1}
\sqrt{n}\left(\widehat{\theta}_n-\theta_\mathbf{T}\right)=-\mathbf{S}^{-1}
\sqrt{n}\mathbf{\Phi}_n(\theta_\mathbf{T})+o_P(1).
\end{equation}
On the other hand, expanding in Taylor series
$$\frac{2n}{\phi^{\prime
\prime}(1)}\widehat{D_\phi}(\theta_0,\widehat\theta_n)=\frac{2n}{\phi^{\prime
\prime}(1)}\int_{\mathbf{I}}
\mathbf{m}(\widehat\theta_n,\mathbf{u})~d\mathbf{C}_{n}(\mathbf{u})$$
in $\widehat\theta_n$ around $\theta_\mathbf{T}$, in connection
with the fact that $\int_{\mathbf{I}}
\mathbf{m}(\theta_\mathbf{T},\mathbf{u})~d\mathbf{C}_{n}(\mathbf{u})=0$,
we get
\begin{equation*}
\frac{2n}{\phi^{\prime
\prime}(1)}\widehat{D_\phi}(\theta_0,\widehat\theta_n)=\frac{2n}{\phi^{\prime
\prime}(1)}\mathbf{\Phi}_n(\theta_\mathbf{T})(\widehat{\theta}_n-\theta_\mathbf{T})-
\frac{n}{\phi^{\prime
\prime}(1)}(\widehat{\theta}_n-\theta_\mathbf{T})^\top\mathbf{\Upsilon}_n(\theta_\mathbf{T})
(\widehat{\theta}_n-\theta_\mathbf{T})+o_P(1).
\end{equation*}
Using (\ref{vv}), (\ref{eqn 1}) and the fact that
$\mathbf{S}=\phi^{''}(1)\mathds{I}_{\theta_\mathbf{T}}$
($\mathds{I}_{\theta_\mathbf{T}}$ denotes the Fisher information
matrix) when $\theta_\mathbf{T}=\theta_0$ to obtain
\begin{equation*}
\frac{2n}{\phi^{\prime
\prime}(1)}\widehat{D_\phi}(\theta_0,\widehat\theta_n)=\frac{1}{\phi^{\prime
\prime}(1)}\sqrt{n}\mathbf{\Phi}_n(\theta_\mathbf{T})\mathds{I}_{\theta_\mathbf{T}}^{-1}\sqrt{n}\mathbf{\Phi}_n(\theta_\mathbf{T}).
\end{equation*}
Finally, use the convergence in (\ref{nom}) and the fact that
$\mathbf{M}=\phi^{''}(1)\mathds{I}_{\theta_\mathbf{T}}$  when
$\theta_\mathbf{T}=\theta_0$, to conclude that
$\frac{2n}{\phi^{\prime
\prime}(1)}\widehat{D_\phi}(\theta_0,\theta_\mathbf{T})$ converges
in distribution to a $\chi^2$ variable with $p$ degrees of freedom
when $\theta_\mathbf{T}=\theta_0$. \hfill$\blacksquare$  \vskip5pt

\noindent (2) Assume that $\theta_\mathbf{T}\neq \theta_0$, using
Taylor expansion again of
$$\widehat{D_\phi}(\theta_\mathbf{T},\theta_0)=\int_{\mathbf{I}}
\mathbf{m}(\widehat\theta_n,\mathbf{u})~d\mathbf{C}_{n}(\mathbf{u})$$
in $\widehat{\theta}_n$ around $\theta_\mathbf{T}$, combined with
the fact that $$\int_{\mathbf{I}} \frac{\partial}{\partial\theta}
\mathbf{m}(\theta_\mathbf{T},\mathbf{u})~d\mathbf{C}_{\theta_\mathbf{T}}(\mathbf{u})=0,$$
we obtain from part (2) of Theorem \ref{theoreme1}
$$\int_{\mathbf{I}} \mathbf{m}(\widehat\theta_n,\mathbf{u})~d\mathbf{C}_{n}(\mathbf{u})=\int_{\mathbf{I}}
\mathbf{m}(\theta_\mathbf{T},\mathbf{u})~d\mathbf{C}_{n}(\mathbf{u})+o_P(n^{-1/2}).$$
Hence,
\begin{eqnarray*}
\lefteqn{\sqrt{n}\left(\widehat{
D_\phi}(\theta_0,\theta_\mathbf{T})-D_\phi(\theta_0,\theta_\mathbf{T})\right)=}\\
&&\sqrt{n}\left(\int_{\mathbf{I}}
\mathbf{m}(\theta_\mathbf{T},\mathbf{u})~d\mathbf{C}_{n}(\mathbf{u})-\int_{\mathbf{I}}
\mathbf{m}(\theta_\mathbf{T},\mathbf{u})~d\mathbf{C}_{\theta_\mathbf{T}}(\mathbf{u})\right)+o_P(1),
\end{eqnarray*}
which under assumption (C.2.iii) by  Theorem 2.1 in
\cite{Ruymgaart_Shorack_Zwet1972} once more, converges to a
centred normal variable with variance given in (\ref{mm}).
\hfill$\blacksquare$  \vskip5pt

\section*{ACKNOWLEDGEMENT}
\small
We would like to thank the
Editor, an Associate editor, and the referees for their
constructive criticism and helpful comments. We are also grateful
to Professors Michel Broniatowski and Paul Deheuvels  for their
helpful discussions and suggestions leading to improvement of this
paper.
\footnotesize
\begin{flushright}
(Revised \today.)
\end{flushright}


\end{document}